\documentclass[10pt]{article}
\usepackage{amssymb}
\usepackage{amsmath,amsthm}
\usepackage{mathtools}
\usepackage{bm}
\usepackage{extarrows}
\usepackage{indentfirst}
\usepackage[english]{babel}
\usepackage{tikz}
\usepackage[numbers]{natbib}
\usepackage[utf8]{inputenc}

\makeatletter

\newcommand{\Rmnum}[1]{\expandafter\@slowromancap\romannumeral #1@}
\makeatother

\def\ex{\mbox{ex}}

\newtheorem{thm}{Theorem}[section]
\newtheorem{cor}[thm]{Corollary}
\newtheorem{lem}[thm]{Lemma}
\newtheorem{prop}[thm]{Proposition}
\newtheorem{conj}[thm]{Conjecture}

\newtheorem{obse}[thm]{Observation}
\newtheorem{claim}{Claim}[section]

\begin{document}
	  \title{Some exact results of the generalized Tur\'an numbers for paths~\thanks{The work was supported by National Natural Science Foundation of China (No. 12071453) and the National Key R and D Program of China(2020YFA0713100).}}
	  \author{Doudou Hei$^1$, Xinmin Hou$^{1,2}$, Boyuan Liu$^1$\\
        \small $^1$School of Mathematics Sciences\\
	  	\small University of Science and Technology of China, Hefei, Anhui 230026, China\\
\small $^{2}$CAS Key Laboratory of Wu Wen-Tsun Mathematics\\
\small University of Science and Technology of China, Hefei, Anhui 230026, China
}
	  \date{}
	  \maketitle
	  	

      \begin{abstract}
For graphs $H$ and $F$ with chromatic number $\chi(F)=k$,
we call $H$ strictly $F$-Tur\'an-good (or $(H, F)$ strictly Tur\'an-good) if the Tur\'an graph $T_{k-1}(n)$ is the unique $F$-free graph on $n$ vertices containing the  largest number of copies of $H$ when $n$ is large enough. Let  $F$ be a graph with chromatic number $\chi(F)\geq 3$ and  a color-critical edge and let $P_\ell$ be a path with $\ell$ vertices. Gerbner and Palmer (2020, arXiv:2006.03756) showed that $(P_3, F)$ is strictly Tur\'an good if $\chi(H)\ge 4$ and they conjectured that (a) this result is true when $\chi(F)=3$, and, moreover, (b) $(P_\ell, K_k)$ is Tur\'an-good for every pair of integers $\ell$ and $k$.
In the present paper, we show that $(H, F)$ is strictly Tur\'an-good when $H$ is a bipartite graph with matching number $\nu(H)=\lfloor \frac{|V(H)|}{2}\rfloor$ and $\chi(F)= 3$, as a corollary, this result confirms the conjecture (a); we also prove that $(P_\ell, F)$ is strictly Tur\'an-good for $2\le\ell\leq 6$ and $\chi(F)\ge 4$, this also confirms the conjecture (b) for $2\le\ell\leq 6$ and $k\ge 4$.

	\end{abstract}

       \section{Introduction}
        Fix a graph $F$, we say that a graph $G$ is $F$-$free$ if it does not contain $F$ as a subgraph. Let $P_\ell$, $C_\ell$ and $K_\ell$ denote a path, cycle and complete graph on $\ell$ vertices, respectively. Fix graphs $H$ and $G$, we denote the number of copies of $H$ in $G$ by $N(H,G)$. We say that an edge $e$ of a graph $F$ is \emph{color-critical} if deleting $e$ from $F$ results in a graph with smaller chromatic number. Let $e(G)$ denote the number of edges of $G$ and write $\nu(G)$ and $\chi(G)$ for the matching number and chromatic number of $G$, respectively. Throughout the paper, let $[k]=\{1, 2, \ldots, k\}$ for a positive integer $k$.

        A fundamental result in extremal graph theory is the Tur\'an Theorem. It states that the Tur\'an graph $T_{k-1}(n)$, which is the complete $(k-1)$-partite graph on $n$ vertices, where each partite class has cardinality $\left\lfloor\frac{n}{k-1}\right\rfloor$ or $\left\lceil\frac{n}{k-1}\right\rceil$, has the largest number of edges among all the $K_k$-free graph. Generally, Tur\'an theory deals with the function $\ex(n, F)$, which is the largest number of edges in $n$-vertex $F$-free graphs.
       We call an $n$-vertex $F$-free graph with $\ex(n, F)$ edges an \emph{extremal graph} for $F$. See, for example, \cite{S97} for a survey.

        For two graphs $H$ and $F$, the {\it generalized Tur\'an number} $\ex(n, H, F)$ is the largest number of copies of $H$ in an $F$-free graph on $n$ vertices, i.e., $$\ex(n, H, F)=\max\{N(H, G) :  G \text{ is an $F$-free graph on $n$ vertices}\}.$$
        Similarly, we call an $n$-vertex $F$-free graph with $\ex(n, H, F)$ copies of $H$ an \emph{extremal graph}.
       After several sporadic results, Alon and Shikhelman~\cite{AS16} studied this problem  systematically. Since then, this problem has attracted many researchers, see e.g.~\cite{EGMS19,GGMV19,GP19,GP20,GS20,MQ20,MN21,QXG21}.


        However, there are not many pairs $(H, F)$ for which the exact values of $\ex(n, H, F)$ were determined.  As pointed by Gerbner and Palmer~\cite{GP20}, there are few $F$-free graphs that are good candidates for being extremal constructions for maximizing copies of $H$,  an exception is the Tur\'an graph, they call  $H$ to be $F$-Tur\'an-good under this situation.
       More precisely, given a $k$-chromatic graph $F$ and a graph $H$ that does not contain $F$ as a
subgraph, we say that $H$ is {\it $F$-Tur\'an-good} (or {\it strictly $F$-Tur\'an-good}) if $\ex(n, H, F) = N (H, T_{k-1}(n))$ (and the Tur\'an graph $T_{k-1}(n)$ is the unique extremal graph) for every $n$ large enough,
we also call $(H, F)$ to be {\it Tur\'an-good} (or {\it strictly Tur\'an-good}) for short.
Here is a list of some Tur\'an-good pairs as we have known so far:
\begin{itemize}
\item[(i)] (Zykov~\cite{Z49}, Erd\H{o}s~\cite{Er62}) $(K_r, K_k)$ is strictly Tur\'an-good for $2\le r<k$;

\item[(ii)] (Simonovits~\cite{S68}, Ma and Qiu~\cite{MQ20}) $(K_r, F)$ is strictly Tur\'an-good, where $F$ is a graph with $\chi(F)>r\ge 2$ and a color-critical edge;

\item[(iii)] (Gy{\H o}ri, Pach and Simonovits~\cite{GPS91}) $(H, K_3)$ is strictly Tur\'an-good, where $H$ is a bipartite graph with matching number $\lfloor \frac{|V(H)|}{2}\rfloor$ (including the path $P_\ell$,  the even cycle $C_{2\ell}$ and the Tur\'an graph $T_2(m)$);

\item[(iv)] (Gy{\H o}ri, Pach and Simonovits~\cite{GPS91}) $(K_{2,t}, K_r)$ are strictly Tur\'an-good for $t=2,3$.

\item[(v)] (Gerbner and Palmer~\cite{GP20})  $(H, K_k)$ is Tur\'an-good for $k\ge k_0$, where $H$ is a complete multipartite graph and $k_0$ is a constant depending on $H$, and they  also conjecture that this result is true for any graph $H$.

\item[(vi)] (Gerbner and Palmer~\cite{GP20}) $(P_4, C_5)$ and $(C_4, C_5)$ are Tur\'an-good, and in general, they~\cite{GP20} conjectured that $(P_k, C_{2\ell+1})$ and $(C_{2k}, C_{2\ell+1})$ are Tur\'an-good. The asymptotic version of this conjecture has been proved by Gerbner et al~\cite{GGMV19}.

\item[(vii)] (Gerbner and Palmer~\cite{GP20}) $(P_4, B_2)$, $(C_4, B_2)$ and $(C_4, F_2)$ are Tur\'an-good, where $B_k$ (resp. $F_k$) is the graph of $k$ triangles all sharing exactly one common edge (resp. one common vertex).

\item[(viii)] (Gerbner and Palmer~\cite{GP20}) $(P_3, F)$ is Tur\'an-good, where $F$ is a graph with chromatic number $\chi(F)=k\ge4$ and a color-critical edge.

\item[(ix)] (Murphy and Nir~\cite{MN21}, Qian et al~\cite{QXG21}) $(P_4, K_{k})$ and $(P_5, K_k)$ are  Tur\'an-good for $k\ge 4$.

\end{itemize}








 Gerbner and Palmer~\cite{GP20} conjectured that the result (viii) is still holds when the chromatic number $\chi(F)=3$ and they also proposed a conjecture that $(P_{\ell}, K_k)$ is Tur\'an-good for every pair of integers $\ell$ and $k$.
\begin{conj}[Gerbner and Palmer~\cite{GP20}]\label{CONJ: conj1}
$(P_3, F)$ is Tur\'an good if $F$ is a graph with chromatic number $\chi(F)\ge 3$ and a color-critical edge.
\end{conj}
\begin{conj}[Gerbner and Palmer~\cite{GP20}]\label{CONJ: conj2}
For every pair of integers $\ell$ and $k$, $(P_\ell, K_k)$ is Tur\'an-good.
\end{conj}

In this article, we first confirm Conjecture~\ref{CONJ: conj1}. In fact, we give a generalized version of (iii) due to Gy{\H o}ri, Pach and Simonovits~\cite{GPS91}, Conjecture~\ref{CONJ: conj1} can be confirmed as a corollary. The following is the first main result.

\begin{thm}\label{THM: main1}
Let $F$ be a  graph with $\chi(F)=3$ and a color-critical edge and let $H$ be a bipartite graph with matching number $\left\lfloor \frac{|V(H)|}{2}\right\rfloor$.  Then $H$ is strictly $F$-Tur\'an good, i.e., $\ex(n, H, F)=N(H, T_2(n))$ for every $n$ large enough. Moreover, the Tur\'an graph $T_2(n)$ is the unique extremal graph for $(H, F)$.
 \end{thm}

So Conjecture~\ref{CONJ: conj1} is a straight-forward corollary of Theorem~\ref{THM: main1} and (viii).

\begin{cor}\label{COR: conj1}
If $F$ is a graph with chromatic number $\chi(F)\ge 3$ and a color-critical edge, then $(P_3, F)$ is Tur\'an-good.
\end{cor}

{\noindent\bf Remarks:} Early this year, Prof. Gerbner told us that he has proved Conjecture~\ref{CONJ: conj1} in~\cite{G21DMGT} by the progressive induction method of Simonovits for generalized Tur\'an problems, in the same paper and several other references provided by Prof. Gerbner, more Tur\'an-good pairs have been proved, we  continue to list them here.
\begin{itemize}
	\item[(x)] (Gerbner~\cite{G21DMGT}) $(M_\ell, F)$ is Tur\'an-good, where $F$ is a graph with a color-critical edge.
	
	\item[(xi)] (Gerbner~\cite{G21DM}) For any positive integers $m$ and $\ell$, $(P_m, C_{2\ell+1})$ and $(C_{2m}, C_{2\ell+1})$ are Tur\'an-good, these results resolved the conjecture (see (vi)) proposed by Gerbner and Palmer~\cite{GP20}; $(P_m, B_t)$ is Tur\'an-good, where $B_t$ is defined as in (vii).
	
	\item[(xii)] (Gerbner~\cite{G21arX,G22arX2}, Gy{\H o}ri, Wang, Woolfson~\cite{GWW21}) $(K_{a,b}, F)$ and $(S_{a,b}, F)$ are Tur\'an-good, where  $F$ is a 3-chromatic graph with a color-critical edge, $K_{a,b}$ and $S_{a,b}$ are a complete bipartite graph and a double star with $|a-b|\le 1$, respectively.
	
\end{itemize}



For Conjecture~\ref{CONJ: conj2}, Murphy and Nir~\cite{MN21}, and Qian et al~\cite{QXG21}) have confirmed this conjecture for $P_4, P_5$ and $k\ge 4$, we continue to confirm this conjecture for $P_6$ and $k\ge 4$ by showing the following a little more generalized result.

\begin{thm}\label{THM: main3}
Let $F$ be a graph with chromatic number $\chi(F)=k\geq 4$ and a color-critical edge. Then the following holds.

(a) If $2\le\ell\le 6$, then $(P_\ell, F)$ is strictly Tur\'an-good, i.e., $\ex(n, P_\ell, F)=N( P_\ell, T_{k-1}(n))$, and $T_{k-1}(n)$ is the unique extremal graph for  $n$  large enough.

(b) There is $k_0$ such that $(P_\ell, F)$ is strictly Tur\'an-good for $\chi(F)=k\ge k_0$.
\end{thm}

The rest of this paper is organized as follows. In Section 2 we give some preliminaries. Next, in Section 3, we prove a technical theorem, which is important in the proofs of Theorems~\ref{THM: main1} and \ref{THM: main3}. We will give the proofs of Theorems~\ref{THM: main1} and \ref{THM: main3} in Sections 4 and 5, respectively. Finally, we briefly resume this work
and give some new lines of research in the Conclusions.

\section{Preliminaries}
In this section we will present some definitions and results needed in the subsequent sections.
Fix a graph $H$ and consider a graph $G$. For each $v\in V(G)$, let $d_G(v,H)$ denote the number of copies of $H$ in $G$ containing the vertex $v$, and let $\delta (G,H)=\min_{x\in V(G)}d_G(x,H)$.
       Let $\overline{H(G)} $ denote the number of different embeddings $\varphi: V(H)\rightarrow V(G)$ such that
       \begin{equation}
              \begin{aligned}
                  &(1)\quad v_1\neq v_2 \Rightarrow \varphi(v_1)\neq \varphi(v_2),\\
                  &(2) \quad v_1v_2\in E(H)\Rightarrow \varphi(v_1)\varphi(v_2)\in E(G)
              \end{aligned}
              \nonumber
       \end{equation}
    for every pair $v_1, v_2\in V(H)$. Evidently, $\overline{H(G)}/{N(H,G)}$ is equal to the number of automorphisms of $H$. Hence, in any class of graphs $\mathcal{G}$, $\overline{H(G)}$ and $N(H,G)$ attain their maximum for the same $G\in \mathcal{G}$.
Similarly, for each $v\in V(G)$, let $\overline{d_G}(v,H)$ denote the number of embeddings of $H$ in $G$ containing the vertex $v$, and let $\overline{\delta}(G,H)=\min_{x\in V(G)}\overline{d_G}(x,H)$.

Given a graph $G$, write $e(G)$ for $|E(G)|$. Let $X$, $Y$ be disjoint subsets of $V(G)$. By $G[X,Y]$ we denote the bipartite subgraph of $G$ consisting of all edges that have one endpoint in $X$ and another in $Y$.
For mutually disjoint subsets $V_1, V_{2}, \ldots, V_{k}\subseteq V(G)$, similarly we define $G[V_1, \ldots, V_{k}]$ to be the $k$-partite subgraph of $G$ consisting of all edges in ${\cup}_{1\leq i< j\leq k}E(G[V_i,V_j])$.
      Write $K(V_1, \ldots, V_k)$ for the complete $k$-partite graph with color classes $V_1, \ldots, V_k$ and write $K_{t_1, \ldots, t_k}$ for a complete $k$-partite graph $K(V_1, \ldots, V_k)$ with $|V_i|=t_i$ for $i\in [k]$.

      An $s$ blow-up of a graph $H$ is the graph obtained by replacing each vertex $v$ of $H$ by an independent set $W_v$ of size $s$, and each edge $uv$ of $H$ by a complete bipartite graph between the corresponding two independent sets $W_u$ and $W_v$.

The following result due to Alon and Shikhelman~\cite{AS16} gave an asymptotical value of $\ex(n, H, F)$ for general graphs $H$ and $F$.
\begin{prop}[\cite{AS16}]\label{PROP: 2.1}
 Let $H$ be a fixed graph with $t$ vertices. Then $ex(n,H,F)$ $=\Omega (n^t)$ if and only if $F$ is not a subgraph of a blow-up of $H$. Otherwise, $\ex(n, H, F)\leq n^{t-\epsilon}$ for some $\epsilon=\epsilon(H, F)>0$.
 \end{prop}

We give a simple observation for the graph with chromatic number $k$ and a color-critical edge.
\begin{obse}\label{OBSE: prop 2.2}
 If $F$ is a graph with $\chi(F)=k\ge3$ and a color-critical edge then $F$ is a subgraph of a complete $k$-partite graph with one class of order one.
\end{obse}


The following two classical results in extremal graph theory will be used.

 \begin{lem}[The stability lemma~\cite{Er66}]\label{LEM: STA-Thm2.2}
  Let $F$ be a graph with  $\chi(F)=k\geq 3$. Then, for every $\varepsilon>0$, there exist $\xi =\xi (F, \varepsilon)>0$ and $n_0=n_0(F, \varepsilon)\in \mathbb{N} $ such that the following holds.  If $G$ is an $F$-free graph on $n\geq n_0$ vertices with $e(G)\geq e(T_{k-1}(n))-\xi n^2$, then there exists a partition of $V(G)=V_1 \cup \ldots \cup V_{k-1}$ such that $\sum_{i=1}^{k-1}e(G[V_i])<\varepsilon n^2/2$.
\end{lem}

  \begin{lem}[Erd\H{o}s-Stone-Simonovits Theorem, \cite{ES66}]\label{LEM: ESS66}
For any graph $H$ with $\chi(H)=r$,
$$\ex(n, H)=\left(1-\frac 1{r-1}\right){n\choose 2}+o(n^2).$$.
\end{lem}

The graph removal lemma given by  Erd\H{o}s, Frankl and R\"odl~\cite{EFR86}, initiated from Ruzsa and Szemer\'edi~\cite{RS76}, also plays an important role in our proofs. An improved version has been given by Fox~\cite{Fox11}. One well-known application of the graph removal lemma is in property testing (one can see~\cite{AS16} for more details if interested).
\begin{lem}[The graph removal lemma~\cite{EFR86}]\label{LEM: removallem}
For each $\epsilon>0$ and graph $H$ on $h$ vertices there is $\delta=\delta(\epsilon,H)>0$ such that every graph on $n$ vertices with at
most $\delta n^h$ copies of $H$ can be made $H$-free by removing at most $\epsilon n^2$ edges.
\end{lem}

Our proof will use the (iii) given by Gy{\H o}ri, Pach and Simonovits~\cite{GPS91}, we restate it here.
\begin{thm}[\cite{GPS91}]\label{THM: 3-TG}
 Let $H$ be a bipartite graph with $m\geq3$ vertices and $\nu(H)=\lfloor \frac{m}{2}\rfloor$. Then, for every $K_3$-free graph $G$ with $n>m$ vertices, $N(H,G)\leq N(H, T_2(n))$, and equality holds if and only if $G\cong T_2(n)$.
 \end{thm}


Given a graph $G$, write $\mu(G)$ for the largest eigenvalue of its adjacency matrix, $\omega(G)$ for its clique
number, and $W_k(G)$ for the number of walks of length $k$ in $G$.
Nikiforov~\cite{Niki02} showed that
\begin{lem}[\cite{Niki02}]\label{LEM: spe-walk}
For every graph $G$ and $r\ge 1$,
${\mu^r(G)}\leq \frac{\omega(G)-1}{\omega(G)}W_{r-1}(G)$.
\end{lem}

The following lemma will be used in the proof of Theorem~\ref{THM: main3} (b).
\begin{lem}[Gerbner and Palmer~\cite{GP20}]\label{LEM:Lemma1.9} For any graph $H$ there are integers $k_0$ and $n_0$ such that if $k\ge k_0$ and $n\ge n_0$, then for any complete $(k-1)$-partite $n$-vertex $K$ we have $N(H, K)\leq N(H, T_{k-1}(n))$, and the equality holds if and only if $G\cong T_{k-1}(n)$.
\end{lem}


\section{T-Extremal Case}
We say a graph $H$ has the {\it weak $k$-T-property} if $N(H, K)\le N (H, T_{k-1}(n))$ for every complete $(k-1)$-partite graph $K=K_{t_1,\ldots,t_{k-1}}$ with $t_1+\ldots+t_{k-1}=n$ and every $n$ large enough, and  the equality holds if and only if $K\cong T_{k-1}(n)$.
An $n$-vertex $F$-free graph $G$ is called {\it T-extremal} if $|e(G)-e(T_{k-1}(n))|=o(n^2)$.

\begin{thm}\label{THM: main-12}
Let $F$ be a  graph with $\chi(F)=k\ge3$ and a color-critical edge and let $H$ be a connected graph with $\chi(H)< k$. Suppose every $F$-free  $n$-vertex graph $G$ with $N(H, G)=\ex(n, H, F)$ is T-extremal.
If $H$ has the weak $k$-T-property, then $G\cong T_{k-1}(n)$.

 \end{thm}
\begin{proof}
Suppose $G$ is an $n$-vertex  $F$-free graph with $N(H, G)=\ex(n,H,F)$ and $n$ is large enough.
Then $N(H,G)\geq N(H, T_{k-1}(n))$. Denote $m:=|V(H)|$.
We may assume an additional condition for $G$ that \ ${\delta} (G, H)\geq {\delta} (T_{k-1}(n), H).$  Indeed, we can
assume $n\ge n_0+a\binom{n_0}{m}$ for some sufficiently large $n_0$, where $a= N(H, K_m)$. If G does not satisfy the property, then there is a vertex $v_n\in V(G)$ such that  ${d_G} (v_n, H)\leq { \delta} (T_{k-1}(n), H)-1$. Set $G_n=G$ and let $G_{n-1}=G-{v_n}$. Then  we have
     \begin{equation}
            \begin{aligned}
                  {N(H, G_{n-1})}&= {N(H,G_n)}-{d_{G_n}}(v_n, H)\\
                   &\geq {N(H, T_{k-1}(n))}- {\delta}(H, T_{k-1}(n))+1\\
                   &\geq {N(H, T_{k-1}(n-1))}+1.
            \end{aligned}
            \nonumber
     \end{equation}
Assume that $G_\ell$ on $\ell$ vertices with $${N(H, G_{\ell})}\ge {N(H, T_{k-1}(\ell))}+n-\ell$$ has been defined for some $\ell\le n-1$. If there exists some vertex $v_\ell\in V(G)$ such that  ${d_{G_\ell}}(v_\ell, H)\leq {\delta}(T_{k-1}(\ell))-1$, let $G_{\ell-1}=G_\ell-{v_\ell}$. Then we get
     \begin{equation}
            \begin{aligned}
                   {N(H,G_{\ell-1})}&={N(H, G_\ell)}-{d_G}(v_\ell, H)\\
                   &\geq {N(H, T_{k-1}(\ell))}+n-\ell- {\delta}(T_{k-1}(\ell), H)+1\\
                   &\geq {N(H, T_{k-1}(\ell-1))}+n-\ell+1;
            \end{aligned}
            \nonumber
     \end{equation}
      otherwise, terminate.
      Let $G_s$ be the graph for which the above iteration terminates. So $G_s$ has exactly $s$ vertices and ${\delta}(G_s, H)\ge {\delta} (T_{k-1}(s), H).$
     If $s< n_0$, then we have
     $$a \binom{n_0}{m}> a \binom{s}{m}\geq {N(H,G_s)} \geq {N(H, T_{k-1}(s))}+n-s \ge n-s> n-n_0\geq a \binom{n_0}{m},$$
     a contradiction. So we have a subgraph $G_s$ of sufficiently large order $s(\geq n_0)$ with ${N(H, G_s)}\ge {N(H, T_{k-1}(s))}+n-s$ and ${\delta}(G_s, H)\geq {\delta}(T_{k-1}(s), H)$. If we can show $G_s$ is a subgraph of some complete $(k-1)$-partite graph $K=K_{t_1,\ldots, t_{k-1}}$ with $t_1+\ldots+t_{k-1}=s$, then we have $N(H, G_s)\le N(H, K)$. If  $H$ has the weak $k$-T-property, then
       $$ N(H, T_{k-1}(s))+n-s\le N(H, G_s)\le N(H, K)\le N(H, T_{k-1}(s)).$$
      So we have $n=s$ and $G_s=G\cong T_{k-1}(n)$.  Therefore, since $s$ is large enough, we can do the same analysis on $G_s$ as $G$. For the sake of writing convenience, in the following proof, we still use $G$ to denote $G_s$ and show that $G\subseteq K_{t_1,\ldots, t_{k-1}}$ with $t_1+\ldots+t_{k-1}=n$ and $|t_i-\frac{n}{k-1}|=o(n)$.

Let $V_1$, \ldots, $V_{k-1}$ be a partition of $V(G)$ such that $e(G[V_1])+ \ldots + e(G[V_{k-1}])$ is minimized. Since $G$ is T-extremal, for every $\varepsilon >0$ (we may choose $\varepsilon$ sufficiently small), choose $\xi=\varepsilon/2$, when $n$ is large enough, we have
        \begin{equation}\label{EQN: e01}
        	\begin{aligned}
        	e(G[V_1])+ \ldots +e(G[V_{k-1}])<\varepsilon n^2/2
        	\end{aligned}
        \end{equation}
        and
        \begin{equation}\label{EQN: e02}
        	\begin{aligned}
        	e(G[V_1, \ldots, V_{k-1}])> e(T_{k-1}(n))- \varepsilon n^2.
        	\end{aligned}
        \end{equation}

\begin{claim}\label{CLM: Claim 01}
There exists some $\theta=\theta(\varepsilon)$ with $\lim_{\varepsilon\rightarrow 0}\theta(\varepsilon)=0$ such that $\left\lvert \left\lvert V_i\right\rvert-\frac{n}{k-1}\right\rvert<\theta n$ for all $i\in[k-1]$.
\end{claim}
\begin{proof}
  Let $p_i:=\frac{\left\lvert V_i \right \rvert}{n}\in [0,1] $.  So we have
        \begin{equation}
              \begin{aligned}
                     e(K(V_1, \ldots, V_{k-1}))&= \sum_{1\le i<j\le k-1}p_ip_jn^2\\
                     &\geq e(G[V_1, \ldots, V_{k-1}])\\
                     &\geq e(T_{k-1}(n))- \varepsilon n^2\\
                     &\geq\frac{k-2}{2(k-1)}n^2- 2\varepsilon n^2.
              \end{aligned}
              \nonumber
       \end{equation}
        Since $\sum_{1\le i<j\le k-1}p_ip_j$ is maximal if and only if $p_1=\ldots=p_{k-1}=\frac{1}{k-1}$, and the maximum value is $\frac{k-2}{2(k-1)}$. By  the continuity of $\sum_{1\le i<j\le k-1}p_ip_j$, there exists $\theta=\theta(\varepsilon)>0$ with $\lim_{\varepsilon\rightarrow 0}\theta(\varepsilon)=0$ such that $\vert p_i- \frac{1}{k-1} \vert<\theta.$
        This proves the claim.
\end{proof}

Let $\beta =2k\sqrt{\varepsilon }$ and $B_i=\{ v\in V_i: \left\lvert  N_G(v)\cap V_i\right\rvert>\beta n \}$ for $i\in[k-1]$. Let $B={\cup}_{i=1}^{k-1}B_i$ and $U_i= V_i\setminus B$.
Because $\beta > 2\sqrt{\varepsilon}$, we have $$\left\lvert B \right\rvert < \frac{2\varepsilon n^2}{\beta n} < \frac{\beta n }{2}.$$

\begin{claim}\label{CLM: Claim 02}
 $B=\emptyset$.
\end{claim}
\begin{proof}
 Since $V_1$, \ldots, $V_{k-1}$ is a partition of $V(G)$ with minimum $e(G[V_1])+ \ldots+ e(G[V_{k-1}])$, we get $$\left\lvert  N_G(v)\cap V_j\right\rvert \ge  \left\lvert  N_G(v)\cap V_i\right\rvert,$$ for any $v\in V_i$, $i, j \in [k-1]$, and $j\neq i$. This together with the definition of $B$ show that for any $v\in B$ and every $i \in [k-1]$,  $\left\lvert  N_G(v)\cap V_i\right\rvert>\beta n.$ Since $U_i= V_i\setminus B$ and $\left\lvert B \right\rvert < \frac{\beta n }{2}$, it follows that  $\left\lvert  N_G(v)\cap U_i\right\rvert > \frac{\beta n }{2}$.

 Suppose $B\neq\emptyset$. Consider an arbitrary vertex $v\in B$. Choose a subset $S_i\subseteq N_G(v)\cap U_i$ with  $|S_i|=\frac{\beta n }{2}$ for $i\in[k-1]$ (this can be done since $\left\lvert  N_G(v)\cap U_i\right\rvert > \frac{\beta n }{2}$). By the inequality~(\ref{EQN: e02})
 $$e(G[S_1, \ldots, S_{k-1}])\ge e\left(T_{k-1}\left((k-1)\beta n/2\right)\right)-\varepsilon n^2=e\left(T_{k-1}\left((k-1)\beta n/2\right)\right)-o\left(\left((k-1)\beta n/2\right)^2\right).$$
By Lemma~\ref{LEM: ESS66}, for $n$ large enough, $G[S_1, \ldots, S_{k-1}]$ contains a copy of the complete $(k-1)$-partite graph $K_{b, \ldots, b}$, where $b=\left\lvert V(F)\right\rvert$. As $S_i\subseteq N_G(v)\cap U_i$, by Observation~\ref{OBSE: prop 2.2},  $G[\{v\}, S_1, \ldots, S_{k-1}]$ contains a copy of $F$, which is a contradiction.
\end{proof}

By Claim~\ref{CLM: Claim 02}, for every $v\in V_i$,   $\left\lvert  N_G(v)\cap V_i\right\rvert\leq\beta n$.

\begin{claim}\label{CLM: Claim03}
 There exists $\zeta=\zeta(\varepsilon)$ with $\lim_{\varepsilon\rightarrow 0}\zeta(\varepsilon)=0$ such that $$\left\lvert N(v)\cap V_j\right\rvert\geq\left(\frac{1}{k-1}-\zeta\right)n$$ for every $v\in V_i$ and $j\neq i.$
 \end{claim}
 \begin{proof} By symmetry, it suffices to prove that the claim is true for every $v\in V_1$. Fix a vertex $v\in V_1$, let $q_j=\frac{\left\lvert N(v)\cap V_j\right\rvert}{n}$.  We will show this claim by estimating $\overline{d_G}(v, H)$.

       First let us estimate the number of $H$-embeddings containing $v$ in $G[V_1, \ldots, V_{k-1}]$. 
Since $H$ is a connected graph, we can find an order of $V(H)$ starting with $v$ such that, after $v$, each vertex has at least one earlier neighbor (e.g we can order the vertices of $V(H)$ by the breath-first search).  
Then we have  $\sum_{j=2}^{k-1} q_jn$ ways to pick the second vertex, and always at most $n-|V_i|\le n-\frac{n}{k-1}+\theta n$ ways to pick a new vertex, where $i$ is some integer in $[k-1]$ and the inequality holds since  $\left\lvert V_i\right\rvert\geq \frac{n}{k-1}-\theta n$ by Claim~\ref{CLM: Claim 01}. Thus the number of this kind of embeddings  is at most $ m \cdot \sum_{j = 2}^{k-1}q_jn\cdot ((\frac{k-2}{k-1}+\theta)n)^{m-2}$. 

 Now, for each embedding of $H$ in $G$ that contains $v$ but it is not in $G[V_1, \ldots, V_{k-1}]$, it must contain some edge in $\cup_{j=1}^{k-1}E(G[V_j])$. By~(\ref{EQN: e01}) and $B=\emptyset$, the number of this kind is at most $2e(H)\cdot\beta n\cdot n^{m-2}+2e(H)\cdot \varepsilon  n^2/2\cdot (m-2)n^{m-3}$. So we have
       \begin{equation}
              \begin{aligned}
                     \overline{d_G}(v,H)\leq & m \cdot \sum_{j = 2}^{k-1}   q_jn\cdot \left(\left(\frac{k-2}{k-1}+\theta\right)n\right)^{m-2}+ \left(2\beta + (m-2)\varepsilon\right) e(H)n^{m-1}.
              \end{aligned}
              \nonumber
       \end{equation}
Since
        \begin{equation}
       		\overline{d_G}(v, H)\geq \overline{\delta} (T_{k-1}(n),H)=  m \cdot\left(\frac{k-2}{k-1}\right)^{m-1}\cdot n^{m-1}+o(n^{m-1}),
       	\nonumber
       \end{equation}
we have
       \begin{equation}
              \begin{aligned}
                     \sum_{j = 2}^{k-1}q_j&\geq\frac{ m\cdot \left(\frac{k-2}{k-1}\right)^{m-1} +o(1)-e(H)(2\beta + (m-2)\varepsilon)}{m\cdot \left(\frac{k-2}{k-1}+\theta \right)^{m-2}}=h(\theta, \beta, \varepsilon).
              \end{aligned}
              \nonumber
       \end{equation}
Note that when $\varepsilon\rightarrow 0$ we have $\beta \rightarrow 0$ and $\theta\rightarrow 0$. So, when $n$ is large enough,
       $$\lim_{\varepsilon \to 0}h(\theta, \beta, \varepsilon)=\frac{k-2}{k-1}. $$
Since  $q_j<\frac{1}{k-1}+\theta$ and $\theta=\theta(\varepsilon)$ is small enough, there exists $\zeta=\zeta(\varepsilon)$ with $\lim_{\varepsilon\rightarrow 0}\zeta(\varepsilon)=0$ such that $q_j\geq \frac{1}{k-1}-\zeta,$ $i\in\{ 2, \ldots, k-1\}$. This proves the claim.
\end{proof}

\begin{claim}\label{CLM: Claim 04}
Every $V_i$ is an independent set in $G$ for $i\in[k-1]$.
\end{claim}
\begin{proof} Suppose to the contrary that say, there exists an edge $e=xy$ in $G[V_1]$. Choose $F_1\subseteq V_1$ with $x, y\in F_1$ and $|F_1|=m$. If we can find  a complete $(k-1)$-partite subgraph $K(F_1, \ldots,  F_{k-1})$ in $G$ such that that $F_i\subseteq V_i$ with $\left\lvert F_i\right\rvert=m$, then, by Observation~\ref{OBSE: prop 2.2}, there exists a copy of $F$ in $G$, which is a contradiction.

To do this, suppose inductively that for some $i\in [k-1]$, we have obtained a complete $i$-partite subgraph  $K(F_1, \ldots,  F_i)$ of $G$. Then  the number of common neighbors of $H_i=\cup _{j=1}^{i}F_j$ in $V_{i+1}$ is at least
       \begin{equation}
       	\begin{aligned}
       		&\sum_{v\in H_i}\left\lvert N_G(v)\cap  V_{i+1}\right\rvert-(\left\lvert H_i\right\rvert-1)\left\lvert V_{i+1}\right\rvert\\
            &\geq \left\lvert H_i\right\rvert\left(\frac{1}{k-1}-\zeta\right)n-(\left\lvert H_i\right\rvert-1)\left(\frac{1}{k-1}+\theta\right)n\\
       		&\geq \left(\frac{1}{k-1}-\left\lvert H_i\right\rvert(\theta+\zeta)\right)n\\
            &\geq m,
       	\end{aligned}
       	\nonumber
       \end{equation}
the last inequality holds since $\zeta, \theta$ are sufficiently small and $n$ is sufficiently large. So  we can find the desired $F_{i+1}$ and the proof of the claim is completed.
\end{proof}

Therefore, $G$ must be a subgraph of $K=K(V_1, \ldots, V_{k-1})$ with $||V_i|-\frac{n}{k-1}|=o(n)$. So $N(H, G)\le N(H, K)$. The proof of the theorem is completed.
\end{proof}

\section{Proof of Conjecture~\ref{CONJ: conj1}}
It is sufficient to prove Theorem~\ref{THM: main1}. A bipartite graph $H$ is said to have the {\it strong T-property} with respect to $F$ if for any $F$-free graph $G$ on $n$ vertices,
$N(H, G)\leq N(H, T_{\chi(F)-1}(n))$
for $n$ large enough, and  equality holds if and only if $G\cong T_{\chi(F)-1}(n)$.
Clearly, if $H$ has the strong T-property with respect to $F$, then $H$ has the weak $\chi(F)$-T-property.
\begin{proof}[Proof of Theorem~\ref{THM: main1}:]    Recall that $F$ is a graph with $\chi(F)=3$ and a color-critical edge and $H$ is a bipartite graph with $m\geq3$ vertices and $\nu(H)=\lfloor \frac{m}{2}\rfloor$.
It is sufficient to show that $H$ has the strong T-property with respect to $F$.
Let $G$ be an $n$-vertex $F$-free graph with the largest number of the copies of $H$. Then $\overline{H(G)}\ge \overline{H(T_{2}(n))}$.
We may assume that $H$ is connected. Indeed, let $H_1, \ldots, H_k$ be the connected components of $H$ and let $|V(H_i)|=m_i$ and set $m_0=0$. Since $\nu(H)=\lfloor \frac{m}{2}\rfloor$,  each component $H_i$ has $\nu(H_i)=\lfloor\frac{m_i}2\rfloor$.
If every $H_i$ has the strong T-property with respect to $F$, then we can embed the components of $H$ into $G$ successively and obtain
     $$\overline{H(G)}=\prod \limits_{i=1}^k \overline{H_i(G(n-\sum_{j< i} m_j))}\leq \prod \limits_{i=1}^k \overline{H_i(T_2(n-\sum_{j< i} m_j))}=\overline{H(T_2(n))},$$
where $G(n-\sum_{j< i} m_j)$ is the subgraph of $G$ on $n-\sum_{j< i} m_j$ vertices. So the equality holds for every $i\in [k]$ and the strong T-property of $H_i$ implies that $G(n-\sum_{j< i} m_j)\cong T_2(n-\sum_{j< i} m_j)$ for $i\in [k]$. Therefore,  $G\cong T_2(n)$.
By Theorem~\ref{THM: 3-TG}, $H$ has the strong T-property with respect to $K_3$. So, by Theorem~\ref{THM: main-12}, it is sufficient to  show that $G$ is T-extremal.

Let $I_k=kK_2$ and $I_k^+=I_k\cup K_1$. By the color-critical edge theorem of Simonovits~\cite{S68} (see (ii) in the introduction),  $K_2$ has the strong T-property with respect to $F$. So $I_k$ and $I_k^+$ have the strong T-property with respect to $F$ too.

\begin{claim}\label{CLM: removeedges}
 $G$ is T-extremal, i.e. $e(G)=e(T_{2}(n))-o(n^2). $
\end{claim}
\begin{proof}
 Since $F$ is a subgraph of a blow-up of $K_3$, by Proposition~\ref{PROP: 2.1},
    $$N(K_3,G)\leq \ex(n, K_3, F)=o(n^3).$$
 By Lemma~\ref{LEM: removallem}, we can get a $K_3$-free graph $G^*$ from $G$ by removing $o(n^2)$ edges.  The number of copies of $H$ intersecting the  $o(n^2)$ removed edges is at most $o(n^2)\cdot O(n^{m-2})=o(n^m)$. Hence
 $$N(H, T_2(n))\leq N(H, G)\leq N(H, G^*)+o(n^m)\leq N(H, T_2(n))+o(n^m),$$
 the last inequality holds because $G^*$ is $K_3$-free and $H$ has the strong T-property with respect to $K_3$ by Theorem~\ref{THM: 3-TG}. So $|\overline{H(G^*)}-\overline{H(T_2(n))}|=o(n^m).$


 Let $a_1b_1, \ldots, a_kb_k\in E(H)$ be a maximum matching of $H$ and let $V(H)=A_0\cup\{a_1, b_1, \ldots, a_k,b_k\}$, where $A_0=\{a_0\}$ if $m$ is odd and $\emptyset$ otherwise. Assume without loss of generality that $A_0\cup\{a_1, \ldots, a_k\}$ and $\{b_1, \ldots, b_k\}$ are the colour classes of $H$.
Since $I_k$ (or $I_k^+$) has the strong T-property, there are $\overline{I_k(G^*)}(\leq\overline{I_k(T_2(n))})$ (or $\overline{I_k^+(G^*)}(\leq\overline{I_k^+(T_2(n))}$) injections $\varphi: A_0\cup\{a_1, \ldots, a_k, b_1, \ldots, b_k\}\rightarrow V(G)$ such that $\varphi(a_i)\varphi(b_i)\in E(G)$ for every $i\in[k]$. Two such injections $\varphi_1$ and $\varphi_2$ are called equivalent if
    \begin{equation}
           \begin{aligned}
               & (1)\   \varphi_1(a_1)=\varphi_2(a_1)\  (\text{ or } \varphi_1(a_0)=\varphi_2(a_0) \text{ if $A_0=\{a_0\}$}), \\
               & (2)\   \{\varphi_1(a_i), \varphi_1(b_i)\}=\{\varphi_2(a_i), \varphi_2(b_i)\}  \quad \forall i\in[k].
           \end{aligned}
           \nonumber
    \end{equation}
So  there are exactly $2^{k-1}$ (or $2^k$ if $A_0=\{a_0\}$) elements in every equivalent class. However, due to the fact that $H$ is connected and $G^*$ is triangle-free, each class contains at most one embedding of $H$ into $G^*$. Thus
    $$ \overline{H(G^*)}\leq 2^{1-k}\overline{I_k(G^*)}\leq 2^{1-k}\overline{I_k(T_2(n))}=\overline{H(T_2(n))}.$$
   $$(\text{or } \overline{H(G^*)}\leq 2^{-k}\overline{I_k^+(G^*)}\leq 2^{-k}\overline{I_k^+(T_2(n))}=\overline{H(T_2(n))} \text{ if $A_0=\{a_0\}$}.)$$
So we have $|\overline{H(G^*)}|=|\overline{H(T_2(n))}|-o(n^m)$.
 We claim that there exists an $i\in[k]$ such that $H-{a_ib_i}$ is still connected. Construct a graph $H^*$ as follows. Let $V(H^*)=\{c_1, c_2, \ldots, c_k\}$ and $c_ic_j\in E(H^*)$  if and only if $E(H[\{a_i, b_i, a_j, b_j\}])\neq \emptyset$.
 If $H^*$ is a tree, assume $c_i$ is a leaf, then $H-{a_ib_i}$ is still connected. Otherwise, there exists a cycle $C$ in $H^*$, assume $c_j\in V(C)$, then $H-{a_jb_j}$ is still connected. Without loss of generality, assume $H'= H-{a_1b_1}$ is still connected.
 Assume $e=xy\in E(G^*)$, let $n(a_1\rightarrow x, b_1\rightarrow y)$ denote the number of different embeddings  $\varphi: V(H)\rightarrow V(G^*)$ such that
    \begin{equation}
           \begin{aligned}
               &(1)\quad \varphi(a_1)=x, \varphi(b_1)=y,\\
               &(2)\quad v_1\neq v_2 \Rightarrow \varphi(v_1)\neq \varphi(v_2),\\
               &(3) \quad v_1v_2\in E(H)\Rightarrow \varphi(v_1)\varphi(v_2)\in E(G^*)
           \end{aligned}
           \nonumber
    \end{equation}
    for every pair $v_1, v_2\in E(H)$. Thus we have

    \begin{equation}
           \begin{aligned}
                  \overline{H(G^*)}&=\sum _{xy\in E(G^*)}\left(n(a_1\rightarrow x, b_1\rightarrow y)+n(a_1\rightarrow y, b_1\rightarrow x)\right)\\
                  &=\sum _{xy\in E(G^*)}\overline{H'(G^*-\{x, y\})}\\
                  &\leq \sum _{xy\in E(G^*)} 2^{2-k}\cdot\overline{I_{k-1}(G^*-\{x, y\})}\left(\text{or } \sum _{xy\in E(G^*)} 2^{1-k}\cdot\overline{I_{k-1}^+(G^*-\{x, y\})}\right)\\
                  &\leq |E(G^*)\vert\cdot 2^{2-k}\cdot\overline{I_{k-1}(T_2(n-2))} \left(\text{or } |E(G^*)\vert\cdot 2^{1-k}\cdot\overline{I_{k-1}^+(T_2(n-2))}\right)\\
                  &= |E(G^*)\vert\cdot 2^{2-k}\prod_{i=1}^{k-1}\overline{I_1(T_2(n-2i))} \left(\text{or }|E(G^*)\vert\cdot 2^{1-k}\prod_{i=1}^{k-1}\overline{I_1(T_2(n-2i))}\cdot(n-2k)\right)\\
                  &= |E(G^*)\vert\cdot 2^{2-k}\prod_{i=1}^{k-1}2e(T_2(n-2i)) \left(\text{or }|E(G^*)\vert\cdot 2^{1-k}\prod_{i=1}^{k-1}2e(T_2(n-2i))\cdot(n-2k)\right),\\
           \end{aligned}
           \nonumber
    \end{equation}
where the case of $A_0=\{a_0\}$ are included in the parentheses.
Combining with
\begin{eqnarray*}
\overline{H(G^*)}= \overline{H(T_2(n))}-o(n^m)&=&2^{1-k}\overline{I_k(T_2(n))}-o(n^m)\\
                 &=& 2^{1-k}\prod_{i=1}^{k}2e(T_2(n-2(i-1)))-o(n^m),
\end{eqnarray*}
or when $A_0=\{a_0\}$
\begin{eqnarray*}
\overline{H(G^*)}= \overline{H(T_2(n))}-o(n^m) &=&2^{-k}\overline{I_k^+(T_2(n))}-o(n^m)\\
                 &=& 2^{-k}\prod_{i=1}^{k}2e(T_2(n-2(i-1)))\cdot(n-2k)-o(n^m),
\end{eqnarray*}
we have $$e(G)=e(G^*)+o(n^2)=e(T_{2}(n))-o(n^2).$$
\end{proof}
This completes the proof of Theorem~\ref{THM: main1}.
\end{proof}

 \section{Proof of Theorem~\ref{THM: main3}}

\begin{proof}[Proof of Theorem~\ref{THM: main3}]
Recall that $F$ has $\chi(F)=k\ge 4$ and a color-critical edge. Let $G$ be an  $n$-vertex $F$-free graph with the largest number of copies of $P_\ell$. So $N(P_\ell, G)\ge N(P_\ell, T_{k-1}(n))$.
We first show that $G$ is T-extremal.
\begin{claim}\label{CLM: Proposition 4.1}
$G$ is T-extremal, i.e. $e(G)=e(T_{k-1}(n))-o(n^2)$.
 \end{claim}
       \begin{proof}
       Since every path corresponds to two walks (one starting from each end-vertex of the path), we have $2N(P_\ell, G)\le W_{\ell-1}(G)$.
       It is well known that $$\frac{W_{\ell-1}(G)}n=\frac{{\bf 1}^TA^{\ell-1}{\bf 1}}{{\bf 1}^T{\bf 1}}\le\mu^{\ell-1}(G),$$
where ${\bf 1}$ is the column vector with all entries being 1 and the last inequality holds because the spectral radius of any Hermitian matrix $M$ is the supremum of the quotient $\frac{{\bf x}^TM{\bf x}}{{\bf x}^T{\bf x}}$, where ${\bf x}$ ranges over $C^n\setminus\{\bf 0\}$.
So we have
       $${\mu}^{\ell-1}(G)n\geq W_{\ell-1}(G)\ge 2N(P_\ell, G)\ge 2N(P_\ell, T_{k-1}(n))=\left(1-\frac{1}{k-1}-o(1)\right)^{\ell-1}n^\ell.$$
By Lemma~\ref{LEM: spe-walk}, ${\mu}^2(G)\leq \frac{\omega(G)-1}{\omega(G)}W_1(G)=2\cdot\frac{\omega(G)-1}{\omega(G)}e(G)$. Hence we have
       $$e(G)\geq \frac{1}{2}\frac{\omega(G)}{\omega(G)-1}{\mu}^2(G)\geq \frac{1}{2}{\mu}^2(G)\geq \frac12\left(1-\frac 1{k-1}-o(1)\right)^2n^2=e(T_{k-1}(n))-o(n^2).$$
\end{proof}

Lemma~\ref{LEM:Lemma1.9} tell us that, for the path $P_\ell$, there is a $k_0$ such that $P_\ell$ has the weak T-property for every $k\ge k_0$.
So if $\chi(H)=k\ge k_0$ then $G\cong T_{k-1}(n)$ by Theorem~\ref{THM: main-12}. This completes the proof of Theorem~\ref{THM: main3} (b).

To prove Theorem~\ref{THM: main3} (a), we will apply Theorem~\ref{THM: main-12} to $F$ and $P_\ell$ with $2\le\ell\le 6$. Clearly, $P_\ell$ is a connected bipartite graph with $\nu(P_\ell)=\lfloor\frac \ell2\rfloor$. The (ii) (Simonovits~\cite{S68}), (viii) (Gerbner and Palmer~\cite{GP20}) and (ix) (Murphy and Nir~\cite{MN21}, Qian et al~\cite{QXG21}) imply that $P_\ell$ has the weak T-property when $\ell\le 5$. To get the result of Theorem~\ref{THM: main3} (a), it is sufficient to show that $P_6$ has the weak T-property.

\begin{claim}\label{CLM: Theorem 1.8}
$P_{6}$ has the weak T-property.
\end{claim}
 \begin{proof} It suffices to prove among all complete $(k-1)$-partite graphs $K=K(V_1, \ldots, V_{k-1})$ on $n$ vertices, the Tur\'an graph $T_{k-1}(n)$ is the unique one with the largest number of copies of $P_6$. Suppose to the contrary that there is a complete $(k-1)$-partite graph $K=K(V_1, \ldots, V_{k-1})$ on $n$ vertices with $N(H, K)\ge N(H, T_{k-1}(n))$ but $K\ncong T_{k-1}(n)$. By Claim~\ref{CLM: Proposition 4.1}, we may assume $\vert |V_i|-\frac{n}{k-1}\vert <\theta n$ for some sufficiently small $\theta>0$.
 Without loss of generality, we assume $\vert V_1\vert >\vert V_2\vert +1$.
 Let $a=\vert V_1\vert-1$, $b=\vert V_2\vert$. Then $a>b$. Let $n_i=\vert V_i \vert$ for $3\le i\le k-1$. If moving a vertex $v$ from $V_1$ to $V_2$, the resulting complete $(k-1)$-partite graph $K'$ has more copies of $P_6$ than in $K$, then we have a contradiction and the claim holds.

 Now let us move a vertex $v$ from $V_1$ to $V_2$ and count the number of embeddings  of $P_6$  destroyed and created after the moving, respectively. Let $P=v_1v_2v_3v_4v_5v_6$ be a path of length 5. There are 6 choices of $v$ to be embedded in the path.
       First, let us count the number of embeddings $\phi_1$ destroyed when $v_1$ is embedded to $v$. Then $\phi_1(v_2)$ must be in $V_2$ and so has $b$ choices. Thus $\phi_1(v_3)\in V_1$ or $V_i$ for some $3\le i\le k-1$.
We count the total number $\sharp \phi_1$ of this kind of destroyed embeddings by dividing them into fifteen cases according to the images of $v_3, v_4$ and $v_5$:

(1) $\phi_1(v_3)\in V_1, \phi_1(v_4)\in V_2$ and $\phi_1(v_5)\in V_1$ or $V_i$ for some $3\le i\le k-1$;

(2) $\phi_1(v_3)\in V_1, \phi_1(v_4)\in V_i$ for some $3\le i\le k-1$ and $\phi_1(v_5)\in V_1, V_2$ or $V_j$ for some $j\not=i$, $3\le j\le k-1$;

(3) $\phi_1(v_3)\in V_i$ for some $3\le i\le k-1$, $\phi_1(v_4)\in V_1$ and $\phi_1(v_5)\in V_2$, $V_i$ or $V_j$ for some $j\not=i$, $3\le j\le k-1$;

(4) $\phi_1(v_3)\in V_i$ for some $3\le i\le k-1$, $\phi_1(v_4)\in V_2$ and $\phi_1(v_5)\in V_1$, $V_i$ or $V_j$ for some $j\not=i$, $3\le j\le k-1$;

(5) $\phi_1(v_3)\in V_i$ for some $3\le i\le k-1$, $\phi_1(v_4)\in V_j$  for some $j\not=i$, $3\le j\le k-1$ and $\phi_1(v_5)\in V_1$, $V_2$, $V_i$ or $V_h$ for some $h\not=i,j$, $3\le h\le k-1$.
Therefore, we have
\begin{eqnarray*}
        \sharp\phi_1&=&ba(b-1)\left[(a-1)(n-3-a)+\sum_{i = 3}^{k-1}n_i(n-4-n_i)\right]\\
        &&+\sum_{i = 3}^{k-1}ban_i\left[(a-1)(n-3-a)+(b-1)(n-3-b)+\sum_{j=3\atop j\neq i}^{k-1}n_j(n-4-n_j)\right]\\
        &&+\sum_{i = 3}^{k-1}bn_ia\left[(b-1)(n-3-b)+(n_i-1)(n-3-n_i)+\sum_{j=3\atop j\neq i}^{k-1}n_j(n-4-n_j)\right]\\
        &&+\sum_{i=3}^{k-1}bn_i(b-1)\left[a(n-4-a)+(n_i-1)(n-3-n_i)+\sum_{j=3\atop j\neq i}^{k-1}n_j(n-4-n_j)\right]\\
        &&+2\sum_{3\le i<j\le k-1}bn_in_j\left[a(n-4-a)+(b-1)(n-3-b)+(n_i-1)(n-3-n_i)\right]\\
        &&+6\sum_{3\le i<j<h\le k-1}bn_in_jn_h(n-4-n_h).
\end{eqnarray*}
Similarly, we can get the number of destroyed embeddings $\sharp\phi_2$ and $\sharp\phi_3$  corresponding to the cases $\phi_2(v_2)=v$ and $\phi_3(v_3)=v$, respectively (the expressions of them are shown in the Appendix).  So the  total number of destroyed embeddings
$$\sharp\phi(a, b, n_3, \ldots, n_{k-1}, n)=2(\sharp\phi_1+ \sharp\phi_2+ \sharp\phi_3).$$
Now we count the number of new created embeddings after the moving of $v$. Similarly, let $\psi_1$ denote a created  embedding with $\psi_1(v_1)=v$. Then $\psi_1(v_2)$ must in $V_1$ and so has $a$ choices. Thus $\psi_1(v_3)\in V_2$ or $V_i$ for some $3\le i\le k-1$. So the number of created embeddings $\psi_1$ is equal to the function by exchanging the variants $a$ and $b$ in $\sharp\phi_1$, i.e. $\sharp\psi_1(a,b,n_3,\ldots,n_{k-1},n)=\sharp\phi_1(b,a, n_3,\ldots,n_{k-1})$. Similarly, we define $\psi_i$ to be a created embedding with $\psi_i(v_i)=v$ for $i=2,3$. We also have  $\sharp\psi_i(a,b,n_3,\ldots,n_{k-1},n)=\sharp\phi_i(b,a, n_3,\ldots,n_{k-1})$ for $i=2,3$. Therefore,  the  total number of created embeddings
\begin{eqnarray*}
\sharp\psi=\sharp\psi(a, b, n_3, \ldots, n_{k-1}, n)&=&2\sum_{i=1}^3\sharp\psi_i(a,b,n_3,\ldots, n_{k-1},n)\\
                                         &=&2\sum_{i=1}^3\sharp\phi_i(b,a,n_3,\ldots, n_{k-1},n)\\
                                         &=&\sharp\phi(b, a, p_3, p_4, \ldots, p_{k-1}, n).
\end{eqnarray*}
        By tedious calculation (or calculated with the MATLAB), we get

\begin{multline*}
   \sharp\phi-\sharp\psi=
   (a-b)\left\{\left[2a^4+7a^3b+13a^2b^2+2a^2n^2+7ab^3+13abn^2+2b^4+2b^2n^2+\sum_{i=3}^{k-1}(a^2+b^2)n_i^2\right.\right.\\
       +\left.\sum_{i=3}^{k-1}\left(4n(a+b)n_i^2+\sum_{j=3\atop j\neq i}^{k-1}\left(8n(a+b)+3ab+6nn_j+4n\sum_{\ell=3\atop \ell\neq i, j}^{k-1}n_\ell\right)n_in_j\right)\right]\\
       -\left[4a^3n+20a^2bn+20ab^2n+4b^3n+\sum_{i = 3}^{k-1}\left(abn_i+4(a^2+ab+b^2)\sum_{j = 3\atop j\neq i}^{k-1}n_j\right)n_i\right.\\
      \left.+\sum_{i=3}^{k-1}\sum_{j=3\atop j\neq i}^{k-1}\left(4(a+b)n_j+4n_i^2+n_in_j+\sum_{\ell=3\atop \ell\neq i,j}^{k-1}4n_\ell^2\right)n_in_j+\left.(3a+3b)\sum_{i = 3}^{k-1}n_i^3\right]+o(n^4)\right\}.
\end{multline*}
Let $Q$ and  $R$  be the expressions in the first and second square brackets, respectively. Since $\vert x-\frac{n}{k-1}\vert <\theta n$ for $x\in\{a,b\}$ and $\vert n_i-\frac{n}{k-1}\vert <\theta n$ for $3\le i\le k-1$,
we have
   $$Q/n^4>(4k^3-26k^2+42k)\left(\frac{1}{k-1}-\theta\right)^3+17\left(\frac{1}{k-1}-\theta\right)^2+(3k^2-17k+55)\left(\frac{1}{k-1}-\theta\right)^4$$
   and $$ R/n^4<(4k^3-23k^2+22k+33)\left(\frac{1}{k-1}+\theta\right)^4+48\left(\frac{1}{k-1}+\theta\right)^3.$$
       Let $f(k,\theta)=Q/n^4-R/n^4$. Then
       \begin{eqnarray*} f(k,0)&=&\frac{4k^3-26k^2+42k}{(k-1)^3}+\frac{17}{(k-1)^2}+\frac{3k^2-17k+55}{(k-1)^4}-\frac{4k^3-23k^2+22k+33}{(k-1)^4}-\frac{48}{(k-1)^3}\\
       &=&\frac{4k^4-24k^3+111k^2-163k+87}{(k-1)^4}.
       \end{eqnarray*}
It is easy to check that $f(k,0)\ge f(4,0)>0$ when $k\geq4$. By the continuity of $f(k,\theta)$ with respect to $\theta$, we have $f(k,\theta)>0$ when $\theta$ is small enough. So $\sharp\phi>\sharp\psi$ when $n$ is large enough, which is a contradiction.
The proof of the claim is completed.
\end{proof}
\end{proof}

\section{Concluding Remarks}
In this article, we show that  $H$ is strictly $F$-Tur\'an good for graph $F$ with $\chi(F)=3$ and a color-critical edge and bipartite graph $H$ with matching number $\left\lfloor \frac{|V(H)|}{2}\right\rfloor$ (Theorem~\ref{THM: main1}). This result confirms Conjecture~\ref{CONJ: conj1} proposed by Gerbner, C. Palmer~\cite{GP20}. But for Conjecture~\ref{CONJ: conj2}, it is far from being resolved. By Theorem~\ref{THM: main-12}, it is sufficient to show that every $P_\ell$ ($\ell\ge 2$) has the weak T-property.   We leave this as an open problem.
It has been shown that $P_2$ (Simonovits~\cite{S68}), $P_3$ (Gerbner and Palmer~\cite{GP20}), $P_4$ (Murphy and Nir~\cite{MN21}), $P_5$ (Qian et al~\cite{QXG21}), and $P_6$ (Theorem~\ref{THM: main3}) have the weak T-property.

\section{Acknowledgment}
We thank Professor D\'aniel Gerbner for providing us a more general version of Theorem~\ref{THM: main-12}.

\begin{appendix}{\centerline{{\bf Appendix}: Expressions of $\sharp\phi_2$ and $\sharp\phi_3$}}
Let $\phi_2$ be  the number of embeddings  destroyed when $v_2$ is embedded to $v$. Then at least one of $\phi_2(v_1)$ and $\phi_2(v_3)$ must be in $V_2$.
We count the total number $\sharp \phi_2$ of this kind embeddings by dividing them into twenty-five cases according to the images of $v_1,v_3, v_4$ and $v_5$:

(1) $\phi_2(v_1)\in V_2, \phi_2(v_3)\in V_2, \phi_2(v_4)\in V_1$ and $\phi_2(v_5)\in V_2$ or $V_i$ for some $3\le i\le k-1$;

 (2) $\phi_2(v_1)\in V_2, \phi_2(v_3)\in V_2, \phi_2(v_4)\in V_i$ for some $3\le i\le k-1$ and $\phi_2(v_5)\in V_1, V_2$ or $V_j$ for some $3\le j\le k-1$ and $j\not=i$;

(3) $\phi_2(v_1)\in V_2, \phi_2(v_3)\in V_i$ for some $3\le i\le k-1$, $\phi_2(v_4)\in V_1$ and $\phi_2(v_5)\in V_2$, $V_i$ or $V_j$ for some $j\not=i$, $3\le j\le k-1$;

(4) $\phi_2(v_1)\in V_2, \phi_2(v_3)\in V_i$ for some $3\le i\le k-1$, $\phi_2(v_4)\in V_2$ and $\phi_2(v_5)\in V_1$, $V_i$ or $V_j$ for some $j\not=i$, $3\le j\le k-1$;

(5) $\phi_2(v_1)\in V_2, \phi_2(v_3)\in V_i$ for some $3\le i\le k-1$, $\phi_2(v_4)\in V_j$  for some $j\not=i$, $3\le j\le k-1$ and $\phi_2(v_5)\in V_1$, $V_2$, $V_i$ or $V_\ell$ for some $\ell\not=i,j$, $3\le \ell\le k-1$;

(6) $\phi_2(v_1)\in V_i$ for some $3\le i\le k-1$, $\phi_2(v_3)\in V_2$, $\phi_2(v_4)\in V_j$  for some $j\not=i$, $3\le j\le k-1$ and $\phi_2(v_5)\in V_1$, $V_2$, $V_i$ or $V_\ell$ for some $\ell\not=i,j$, $3\le \ell\le k-1$ (the number of the destroyed embeddings is the same as in the case (5));

(7) $\phi_2(v_1)\in V_i$ for some $3\le i\le k-1$, $\phi_2(v_3)\in V_2$, $\phi_2(v_4)\in V_1$ and $\phi_2(v_5)\in V_2$, $V_i$ or $V_j$ for some $j\not=i$, $3\le j\le k-1$;

(8) $\phi_2(v_1)\in V_i$ for some $3\le i\le k-1$, $\phi_2(v_3)\in V_2$, $\phi_2(v_4)\in V_i$ and $\phi_2(v_5)\in V_1$, $V_2$ or $V_j$ for some $j\not=i$, $3\le j\le k-1$.
Therefore, we have
\begin{eqnarray*}
 \sharp\phi_2&=& b(b-1)a\left[(b-2)(n-2-b)+\sum_{i=3}^{k-1}n_i(n-4-n_i)\right]\\
             &&+\sum_{i=3}^{k-1}b(b-1)n_i\left[a(n-4-a)+(b-2)(n-2-b)+\sum_{j=3\atop j\not=i}^{k-1}n_j(n-4-n_j)\right]\\
             &&+\sum_{i=3}^{k-1}bn_ia\left[(b-1)(n-3-b)+(n_i-1)(n-3-n_i)+\sum_{j=3 \atop j\not=i}^{k-1}n_j(n-4-n_j)\right]\\
             &&+\sum_{i=3}^{k-1}bn_i(b-1)\left[a(n-4-a)+(n_i-1)(n-3-n_i)+\sum_{j=3\atop j\not=i}^{k-1}n_j(n-4-n_j)\right]\\
             &&+2\cdot 2\sum_{3\le i<j\le k-1}bn_in_j\left[a(n-4-a)+(b-1)(n-3-b)+(n_i-1)(n-3-n_i)\right]\\
             &&+2\cdot 6\sum_{3\le i<j<\ell\le k-1}bn_in_jn_\ell(n-4-n_\ell)\\
        &&+\sum_{i=3}^{k-1}n_iba\left[(b-1)(n-3-b)+(n_i-1)(n-3-n_i)+\sum_{j=3\atop j\neq i}^{k-1}n_j(n-4-n_j)\right]\\
        &&+\sum_{i=3}^{k-1}n_ib(n_i-1)\left[a(n-4-a)+(b-1)(n-3-b)+\sum_{j=3\atop j\neq i}^{k-1}n_j(n-4-n_j)\right].
 \end{eqnarray*}
Similarly, if $\phi_3(v_3)=v$, then at least one of $\phi_3(v_2), \phi_3(v_4)$ must be in $V_2$. We also can count the total number $\sharp\phi_3$ of this kind of destroyed embeddings by dividing them into twenty-five cases according to the images of $v_2, v_4, v_5$ and $v_1$:

(1) $\phi_3(v_2)\in V_2, \phi_3(v_4)\in V_2, \phi_2(v_5)\in V_1$ and $\phi_2(v_6)\in V_2$ or $V_i$ for some $3\le i\le k-1$;

 (2) $\phi_3(v_2)\in V_2, \phi_3(v_4)\in V_2, \phi_3(v_5)\in V_i$ for some $3\le i\le k-1$ and $\phi_3(v_6)\in V_1, V_2$ or $V_j$ for some $3\le j\le k-1$ and $j\not=i$;

(3) $\phi_3(v_2)\in V_2, \phi_3(v_4)\in V_i$ for some $3\le i\le k-1$, $\phi_3(v_5)\in V_1$ and $\phi_3(v_6)\in V_2$, $V_i$ or $V_j$ for some $j\not=i$, $3\le j\le k-1$;

(4) $\phi_3(v_2)\in V_2, \phi_3(v_4)\in V_i$ for some $3\le i\le k-1$, $\phi_3(v_5)\in V_2$ and $\phi_3(v_6)\in V_1$, $V_i$ or $V_j$ for some $j\not=i$, $3\le j\le k-1$;

(5) $\phi_3(v_2)\in V_2, \phi_3(v_4)\in V_i$ for some $3\le i\le k-1$, $\phi_3(v_5)\in V_j$  for some $j\not=i$, $3\le j\le k-1$ and $\phi_3(v_6)\in V_1$, $V_2$, $V_i$ or $V_\ell$ for some $\ell\not=i,j$, $3\le \ell\le k-1$;

(6) $\phi_3(v_2)\in V_i$ for some $3\le i\le k-1$, $\phi_3(v_4)\in V_2$, $\phi_3(v_5)\in V_j$  for some $j\not=i$, $3\le j\le k-1$ and $\phi_3(v_6)\in V_1$, $V_2$, $V_i$ or $V_\ell$ for some $\ell\not=i,j$, $3\le \ell\le k-1$ (the number of the destroyed embeddings is the same as in the case (5));

(7) $\phi_3(v_2)\in V_i$ for some $3\le i\le k-1$, $\phi_3(v_4)\in V_2$, $\phi_3(v_5)\in V_1$ and $\phi_3(v_6)\in V_2$, $V_i$ or $V_j$ for some $j\not=i$, $3\le j\le k-1$;

(8) $\phi_3(v_2)\in V_i$ for some $3\le i\le k-1$, $\phi_3(v_4)\in V_2$, $\phi_3(v_5)\in V_i$ and $\phi_3(v_6)\in V_1$, $V_2$ or $V_j$ for some $j\not=i$, $3\le j\le k-1$.
\begin{eqnarray*}
\sharp\phi_3&=&b(b-1)a\left[(b-2)(n-2-b)+ \sum_{i=3}^{k-1}n_i(n-3-b)\right]\\
        &&+\sum_{i=3}^{k-1}b(b-1)n_i\left[a(n-3-b)+(b-2)(n-2-b)+\sum_{j=3\atop j\neq i}^{k-1}n_j(n-3-b)\right]\\
        &&+\sum_{i=3}^{k-1}bn_ia\left[(b-1)(n-3-b)+(n_i-1)(n-4-b)+\sum_{j=3 \atop j\neq i}^{k-1}n_j(n-4-b)\right]\\
        &&+\sum_{i=3}^{k-1}bn_i(b-1)\left[a(n-3-b)+(n_i-1)(n-3-b)+\sum_{j=3 \atop j\neq i}^{k-1}n_j(n-3-b)\right]\\
        &&+2\cdot 2\sum_{3\le i<j\le k-1}bn_in_j\left[a(n-4-b)+(b-1)(n-3-b)+(n_i-1)(n-4-b)\right]\\
        &&+2\cdot 6\sum_{3\le i<j<\ell\le k-1}bn_in_jn_\ell(n-4-b)\\
        &&+\sum_{i=3}^{k-1}n_iba\left[(b-1)(n-4-n_i)+(n_i-1)(n-3-n_i)+\sum_{j=3 \atop j\neq i}^{k-1}n_j(n-4-n_i)\right]\\
        &&+\sum_{i=3}^{k-1}n_ib(n_i-1)\left[a(n-3-n_i)+(b-1)(n-3-n_i)+\sum_{j=3 \atop j\neq i}^{k-1}n_j(n-3-n_i)\right].
       \end{eqnarray*}

\end{appendix}

\end{document}